\documentclass[12pt]{amsart} 

\usepackage{microtype,times}
\usepackage{amsmath}
\usepackage{amsthm, amssymb,bm}
\usepackage{graphicx}
\usepackage{color}
\usepackage{todonotes}
\usepackage{listings}
\usepackage{tikz}
\usetikzlibrary{arrows,decorations.pathreplacing}

\usepackage{soul}
\usepackage{ulem}
\usepackage{nicefrac}

\usepackage[capitalize,noabbrev]{cleveref}

\usepackage[backend=bibtex, doi=true, isbn=false, url=false, giveninits=true, style=alphabetic]{biblatex}

\newcommand{\bC}{\mathbb{C}}

\newcommand{\bQ}{\mathbb{Q}}

\newcommand{\bZ}{\mathbb{Z}}

\newcommand{\bN}{\mathbb{N}}

\def\({ \left( }
\def\){ \right)}

\newcommand{\Sym}[1]{\mathfrak{S}_{#1}}
\DeclareMathOperator{\Ch}{Ch}
\newcommand{\ChSpin}{\Ch^{\mathrm{spin}}}

\newcommand{\GA}{{G}}
\newcommand{\GB}{{H}}
\newcommand{\HA}{{I}}
\newcommand{\HB}{{J}}

\DeclareMathOperator{\Cat}{Cat}

\theoremstyle{plain}
\newtheorem{thm}{Theorem}[section]
\newtheorem{prop}[thm]{Proposition}
\newtheorem{lem}[thm]{Lemma}

\newtheorem{cor}[thm]{Corollary}

\theoremstyle{definition}

\theoremstyle{remark}
\newtheorem{remark}[thm]{Remark}
\newtheorem{example}[thm]{Example}
\newtheorem{problem}[thm]{Problem}

\crefname{thm}{Theorem}{Theorems}

\numberwithin{equation}{section}

\title[Characters of almost square shape]{Symmetric group characters \\ of 
almost square shape}
\author{\textsc{Sho Matsumoto}
and \textsc{Piotr \'Sniady}
}

\begin{document}

\begin{abstract}
    We give closed product formulas for the irreducible characters of the
    symmetric groups related to rectangular \emph{`almost square'} Young diagrams
    $p\times(p+\delta)$ for a fixed value of an integer~$\delta$ and an arbitrary integer $p$.
\end{abstract}



\maketitle

\lstset{language=Python}
\lstset{basicstyle=\footnotesize\ttfamily,breaklines=true}

\section{Introduction}

\subsection{The main result}
Let $\pi$ be a partition of an integer $k$ and let $\lambda$ be a partition of an integer~$n$.
Let $\Ch_\pi(\lambda)$ denote \emph{the normalized character} of the symmetric group $\Sym{n}$ defined by 
\begin{equation} \label{def:normalchar}
\Ch_\pi(\lambda)= \begin{cases} n^{\downarrow k} \cdot \frac{\chi^\lambda_{\pi\cup (1^{n-k})}}{f^\lambda},
    & \text{if } n\geq k, \\
    0 & \text{otherwise,}
    \end{cases}
\end{equation}
where $\chi^\lambda_\mu$ is the usual character of the irreducible representation of the symmetric group $\Sym{n}$
associated with $\lambda$, evaluated on the conjugacy class associated with $\mu$.
This choice of the normalization is quite natural, in particular in the context of the asymptotic representation theory, 
see for example \cite{Biane2003,IvanovKerov}.

\begin{figure}[t]
    \centering
    \begin{tikzpicture}
        \draw(0,0) grid (6,5);
        \draw[ultra thick] (0,0) rectangle (6,5); 
        \draw [gray,<->] (0,5.5) -- (6,5.5) node [midway,outer sep=5pt,fill=white] {\textcolor{black}{$q$}};
        \draw [gray,<->] (6.5,0) -- (6.5,5) node [midway,outer sep=5pt,fill=white] {\textcolor{black}{$p$}};
    \end{tikzpicture}

\caption{Rectangular Young diagram $p\times q$.}
\label{fig:rectangular}
\end{figure}
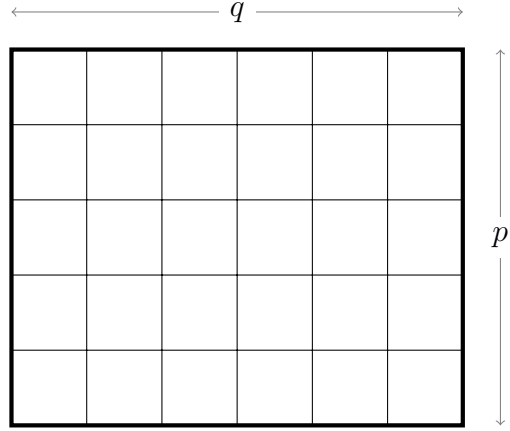

In this note we will concentrate on the case when $\pi=(k)$ consists of a
single part, i.e., on the characters evaluated on a single cycle (augmented by
a necessary number of fixpoints). 
Also, we will concentrate on the special case when 
\[\lambda=p\times q= (\underbrace{q,\dots,q}_{\text{$p$ times}})\] 
is a rectangular Young diagram, see \cref{fig:rectangular}. 
We will give 
closed product formulas for such characters in the \emph{almost square} setting when
$q-p$ is a fixed integer and $p$ is arbitrary.
The exact form of the formula depends on the parity of the length $k$ of the cycle,
as well as on the parity of the difference $q-p$, 
so altogether there are four distinct formulas for such characters. As a teaser, we start with the case when $k=2j-1$ is odd while $q-p=2d$ is even.

\begin{thm} 
\label{thm:main-odd-even}
Let $j$ be a positive integer.
Then
\begin{multline}
\label{eq:main-A}
\Ch_{2j-1}\big( (e-d) \times (e+d) \big)  \\
= (-1)^{j-1} \Cat (j-1) \sum_{k=0}^j f_k(j)
\left( \prod_{r=0}^{k-1} (d^2-r^2) \right)  
\left( \prod_{r=k}^{j-1} (e^2-r^2) \right),
\end{multline}
where
\[ \Cat (j-1) =\frac{(2j-2)!}{(j-1)!\; j!} \]
is the Catalan number; furthermore
$f_0(j)=1$ and 
\begin{equation}
\label{eq:f_k}
f_k(j)= (-1)^k \frac{j^{\downarrow k} \, (2j-1)^{\uparrow \uparrow k}}{k! \, (2k-1)!!} 
\end{equation}
for $1 \le k \le j$.
\end{thm}

Above we used the \emph{double rising factorial} $a^{\uparrow \uparrow k}$ 
(which is somewhat analogous to the \emph{double factorial} $a!!$ in which the factors form an arithmetic progression with the step $2$) 
which is defined by
\[
a^{\uparrow \uparrow k} =\prod_{r=0}^{k-1} (a+2r)  = 2^k \left( \frac{a}{2} \right)^{\uparrow k}
\]
for a complex number $a$ and a positive integer $k$, and by $a^{\uparrow \uparrow 0}:=1$.
Thus, $f_k(j)\in\bQ[j]$ is a polynomial in the variable $j$ of degree $2k$.

Note that in the above result there is \emph{no} assumption that $e-d$ and $e+d$ 
are non-negative integers; in fact $e$ and $d$ can be arbitrary complex numbers.
The reader may feel uneasy about the case when $(e-d)\times (e+d)$ does not make
sense as a Young diagram; later on in \cref{cor:poly-degree}
we will explain why in this case
the left-hand side of \eqref{eq:main-A} still makes sense.

\subsection{The product formula}

In the special case when $|d|$ is a small integer, the formula
\eqref{eq:main-A} takes a simpler form because each summand 
on the right-hand side which corresponds to $k$ such that $k>|d|$ 
is equal to zero; in this way the sum can be taken over 
$k\in\{1,\dots,j \wedge |d|\}\textbf{}$. This observation is especially convenient in the
aforementioned \emph{almost square} setting when we consider the character
corresponding to a rectangular Young diagram $\lambda=p\times q$ in the setup where
$q-p=2d$ is a fixed even integer and 
$p$ is arbitrary. In particular, we get the following closed product form for the character. 

\begin{cor}
\label{cor:product}
Let $j,p,q$ be positive integers; we denote by $n=pq=|p\times q|$ the number of the boxes of the corresponding Young diagram.
Suppose that $q-p$ is an even integer which we denote by $2d:=q-p$.
Then
\begin{equation}
\label{eq:cor12}
\Ch_{2j-1}\big( p \times q \big) 
= (-1)^{j-1} \Cat (j-1)\ \GA_d(j,n)  \prod_{r=0}^{j-|d|-1} \big(n -r(r+2|d|) \big),
\end{equation}
where 
\[
\GA_{d}(j,n)= \sum_{k=0}^{|d|} f_k(j) \left( \prod_{r=0}^{k-1} (d^2-r^2) \right)  
\left( \prod_{r=k}^{|d|-1} (n+d^2-r^2) \right)
\]
with  $f_k(j)$ given, as before, by \eqref{eq:f_k}.
\end{cor}

We can see that, with $d$ fixed,  $\GA_d(j,n)\in\bQ[j,n]$ is a polynomial in the variables $j,n$ of the total 
degree $2\ |d|$ if we declare that the degrees of the variables $j$ and $n$ are given by $\deg j=1$ and $\deg n=2$. 
In fact, $\GA_d(j,n)\in\bZ[j,n]$ is a polynomial with \emph{integer} coefficients, see \cref{prop:integer}.

\begin{example}
\begin{align*}
\GA_0(j,n) &= 1,\\
\GA_{\pm 1}(j,n)&= n+1-j(2j-1), \\
\GA_{\pm 2}(j,n)&= (n+4)(n+3)-4j(2j-1)(n+3) +2j(j-1)(2j-1)(2j+1).
\end{align*}
\end{example}

\subsection{Convention for products}

\cref{cor:product} remains valid in the case when $j\leq |d|$, however in this case
the product on the right hand side of \eqref{eq:cor12} should be understood using the following non-standard convention:
\[    \prod_{r=0}^{l} a_r =
\begin{cases} 
    a_0 \cdots a_l                          & \text{if } l\geq 0, \\
    1                                       & \text{if } l=-1,     \\
    \frac{1}{a_{l+1} a_{l+2} \cdots a_{-2} a_{-1}}     & \text{if } l\leq  -2.
\end{cases}
\]
This convention was chosen in such a way that the identity
\[ \prod_{r=0}^{l+1} a_r = \left[ \prod_{r=0}^{l} a_r \right] \cdot  a_{l+1} \]
holds for any (positive or negative) integer $l$.

\subsection{Collection of results}
Below we present a collection of the results which cover the remaining choices for the parity for the length of the cycle and the difference $q-p$ between the rectangle sides.

\subsubsection{The length of the cycle is odd, the difference of the rectangle sides is odd}

The following result is a counterpart of \cref{thm:main-odd-even} 
which is particularly useful 
for a rectangular Young diagram $\lambda=p\times q$ for which $q-p$ is an odd integer.

\begin{thm}
\label{thm:odd}
Let $j$ be a positive integer.
Then
\begin{multline*}
\Ch_{2j-1}\big( (e-d) \times (e+d) \big) \\
= (-1)^{j-1} \Cat (j-1) \sum_{k=0}^j f_k(j)
\left[ \prod_{r=0}^{k-1} \big(d^2-(r+\tfrac{1}{2})^2 \big) \right]  
\left[ \prod_{r=k}^{j-1} \big(e^2-(r+\tfrac{1}{2})^2 \big)  \right].
\end{multline*}
\end{thm}

\begin{cor}
\label{cor:odd}
Let $j,p, q$ be positive integers and 
set $n=pq$.
Suppose that $q-p$ is an odd integer $2d:=q-p$, where $d\in\left\{\pm \tfrac{1}{2}, 
\pm \tfrac{3}{2}, \dots\right\}$.
Then
\[
\Ch_{2j-1}\big( p\times q \big) 
= (-1)^{j-1} \Cat (j-1)\ \GB_d(j,n)  \prod_{r=0}^{j-|d|-\tfrac{1}{2}} \big(n -r(r+2|d|) \big),
\]
where 
\begin{multline*}
\GB_{d}(j,n)= \sum_{k=0}^{|d|-\tfrac{1}{2}} 
(-1)^k \frac{j^{\downarrow k} \, (2j-1)^{\uparrow \uparrow k}}{k! \, (2k-1)!!} 
  \prod_{r=0}^{k-1} \big(d^2- (r+\tfrac{1}{2})^2 \big)   \\
\times \prod_{r=k}^{|d|-\tfrac{3}{2}} \big(n+d^2- (r+\tfrac{1}{2})^2 \big)  .
\end{multline*}
\end{cor}

\begin{example}
\begin{align*}
\GB_{\pm \tfrac{1}{2}}(j,n)&= 1, \\
\GB_{\pm \tfrac{3}{2}}(j,n)&= n+2 -2j(2j-1), \\
\GB_{\pm \tfrac{5}{2}}(j,n)&= (n+6)(n+4)-6j(2j-1) (n+4)+4j(j-1)(2j-1)(2j+1).
\end{align*}
\end{example}

\subsubsection{The length of the cycle is even, the difference of the rectangle sides is even}

The following result is a direct counterpart of \cref{thm:main-odd-even} 
for the even cycle.

\begin{thm}
\label{thm:even1}
Let $j$ be a positive integer.
Then
\begin{multline*}
\Ch_{2j}\big( (e-d) \times (e+d) \big) \\
= (-1)^{j-1} \binom{2j-1}{j} 
\sum_{k=0}^{j} g_k(j)
2d\left( \prod_{r=1}^{k} (d^2-r^2) \right)  
\left( \prod_{r=k+1}^{j} (e^2-r^2) \right),
\end{multline*}
where
$g_0(j)=1$ and 
\begin{equation*}
g_k(j)= (-1)^{k} \frac{j^{\downarrow k} \, (2j+1)^{\uparrow \uparrow k}}{k! \, (2k+1)!!} 
\end{equation*}
for $1 \le k \le j$.
\end{thm}

\begin{cor}
\label{cor:HA}
Let $j,p,q$ be positive integers and set $n=pq$.
Suppose that $q-p$ is an even integer $2d=q-p$, where $d\in \{ 0, \pm 1, \pm 2, \dots \}$.
Then
\[
\Ch_{2j}\big( p \times q \big) 
= (-1)^{j-1} \binom{2j}{j}   \HA_d(j,n)  \prod_{r=0}^{j-|d|} \big(n -r(r+2|d|) \big),
\]
where  $\HA_0(j,n)=0$ and, for $d \in \{ \pm 1, \pm 2, \dots\}$
\[
\HA_{d}(j,n)=  d\sum_{k=0}^{|d| -1}  
(-1)^{k} \frac{j^{\downarrow k} \, (2j+1)^{\uparrow \uparrow k}}{k! \, (2k+1)!!}   
\left( \prod_{r=1}^{k} (d^2-r^2) \right)  
\left( \prod_{r=k+1}^{|d|-1} (n+d^2-r^2) \right).
\]
\end{cor}

With $d$ fixed,  
$\HA_d(j,n)\in\bQ[j,n]$ is a polynomial in the variables $j,n$ of the total 
degree $2|d|-1$ if we give $\deg j=1$ and $\deg n=2$.

\begin{example}
\begin{align*}
\HA_{0}(j,n) &= 0, \\
\HA_{1}(j,n) &= 1, \\
\HA_{2}(j,n) &= 2(n+3)-2j(2j+1) =2 \big( n-(j-1)(2j+3)\big), \\
\HA_{3}(j,n) &= 3(n+8)(n+5) -8j(2j+1)(n+5)+4j(j-1)(2j+1)(2j+3).
\end{align*} 
For negative integers $d$, we have $\HA_{d}(j,n)=-\HA_{-d}(j,n)$.
\end{example}

\subsubsection{The length of the cycle is even, the difference of the rectangle sides is odd}

The following result is a direct counterpart of \cref{thm:odd} for an even cycle.

\begin{thm}
\label{thm:even2}
Let $j$ be a positive integer.
Then
\begin{multline*}
\Ch_{2j}\big( (e-d) \times (e+d) \big)  \\
= (-1)^{j-1} \binom{2j-1}{j} 
\sum_{k=0}^{j} g_k(j)
2d\left( \prod_{r=1}^{k} \big( d^2-(r-\tfrac{1}{2})^2 \big ) \right)  
\left( \prod_{r=k+1}^{j} \big( e^2-(r-\tfrac{1}{2})^2 \big )\right).
\end{multline*}
\end{thm}

\begin{cor}
\label{cor:HB}
Let $j,p,q$ be positive integers and set $n=pq$. 
Suppose that $q-p$ is an odd integer $2d:=q-p$, where $d\in \{\pm \tfrac{1}{2}, \pm \tfrac{3}{2}, \dots\}$. 
Then 
\[
\Ch_{2j}\big( p \times q \big) 
= (-1)^{j-1} \binom{2j-1}{j} \HB_d(j,n)  \prod_{r=0}^{j-|d|-\tfrac{1}{2} } \big(n -r(r+2|d|) \big),
\]
where 
{\small
\begin{multline*}
\HB_{d}(j,n)= 2d\sum_{k=0}^{|d|-\tfrac{1}{2}} 
(-1)^{k} \frac{j^{\downarrow k} \, (2j+1)^{\uparrow \uparrow k}}{k! \, (2k+1)!!}   
 \left( \prod_{r=1}^{k} \big( d^2- (r-\tfrac{1}{2})^2 \big) \right)  \\
 \times 
\left( \prod_{r=k+1}^{|d|-\tfrac{1}{2}}\big(n+ d^2- (r-\tfrac{1}{2})^2 \big) \right).
\end{multline*}
}
\end{cor}

\begin{example}
\begin{align*}
\HB_{\tfrac{1}{2}} (j,n) &= 1, \\
\HB_{ \tfrac{3}{2}}(j,n) &= 3(n+2) -2j(2j+1) =3n- 2(j-1)(2j+3), \\
\HB_{\tfrac{5}{2}}(j,n) &= 5(n+6)(n+4)-10j(2j+1) (n+4)+4j(j-1)(2j+1)(2j+3).
\end{align*}
For negative half integers $d$, we have $\HB_{d}(j,n)=-\HB_{-d}(j,n)$.
\end{example}

\subsection{Vanishing of some special characters}

The special case of \cref{cor:odd} for $d=\tfrac{3}{2}$ and $p=2j-2$, $q=2j+1$ gives rise to the following somewhat
surprising corollary for which we failed to find an alternative simple proof.
\begin{cor}
    For each integer $j\geq 2$ the irreducible character related to the
rectangular diagram $(2j-2) \times (2j+1)$ vanishes on the cycle of length
$2j-1$, i.e.
    \[ \chi^{(2j-2) \times (2j+1)}_{2j-1,\ 1^{(2j-3)(2j-1)}} =0. \]
\end{cor}

\subsection{The link with spin characters related to the staircase strict partition}
One of the motivations for the current paper was 
the recent progress related to 
the \emph{spin characters of the symmetric groups}
\cite{Schur11}.
On one hand, De Stavola \cite[Proposition 4.18, page 91]{DeStavola2017} 
gave an explicit formula 
for the spin character related to the \emph{staircase strict partition} 
\begin{equation*}
\Delta_p=(p,p-1,p-2,\dots, 2, 1)
\end{equation*}
which has the property that its \emph{double}
\begin{equation*}
D(\Delta_p)= (\underbrace{p+1,p+1,\dots,p+1}_{\text{$p$ times}}) = p \times (p+1)
\end{equation*}
is a rectangular Young diagram which is almost square.
On the other hand, in our recent paper \cite{Matsumoto2020} 
we found an identity which gives a link between the spin characters
and their usual (linear) counterparts
\begin{equation*}
2\ChSpin_{2j-1}(\xi) =
\Ch_{2j-1}(D(\xi))
\end{equation*}
which holds true for any strict partition $\xi$.

By combining these two results one gets a closed product formula 
for the linear character $\Ch_{2j-1}\big(p\times (p+1) \big)$
corresponding to a rectangular Young diagram which is almost square;
this closed formula coincides with the special case of
\cref{cor:odd} for $d=\tfrac{1}{2}$
(in fact the formula in the original paper of De Stavola
has an incorrect sign). 
In his proof, De Stavola employed some computations in Maple;
by turning the argument around our \cref{cor:odd} gives a 
purely algebraic proof of his product formula 
for $\ChSpin_{2j-1}(\Delta_p)$.

\begin{problem}
It would be interesting to extend this result to the 
\emph{almost} staircase strict partition of the form
\[ \Delta_{p,d}=(\underbrace{p,\dots,p}_{\text{$d+1$ times}},p-1,p-2,\dots, 2, 1) \]
and the corresponding spin character $\ChSpin_{2j-1}(\Delta_{p,d})$ in the setting when
$d\geq 0$ is a small integer and $p$ varies.
The above question is equivalent to finding closed formulas for the 
linear characters $\Ch_{2j-1}\big( D(\Delta_{p,d}) \big)$; 
note that the double $D(\Delta_{p,d})$ can be viewed as a large rectangular Young diagram
\mbox{$(p+d) \times (p+d+1)$} of almost square shape 
from which another, small rectangle $d\times (d+1)$ of almost square shape was removed.
\end{problem}

\subsection{Sketch of the proof}

We will start in \cref{sec:toolbox} by
collecting some formulas for the irreducible characters related to rectangular shapes. 
Our strategy towards the proof of \cref{thm:main-odd-even} is threefold. 

Firstly, we will fix the value of an integer $j\geq 1$ and we shall 
investigate the function
\[ (d,e) \mapsto \Ch_{2j-1}\big( (e -d) \times (e+d) \big).\]
We will show that it is a polynomial in the variables $d,e$ of the total degree $2j$.

Secondly, we will show that this polynomial is of the form
\[  \Ch_{2j-1}\big( (e -d) \times (e+d) \big) =
\sum_{k=0}^j c_k(j) \left(\prod_{r=0}^{k-1} (d^2-r^2) \right) \left(\prod_{r=k}^{j-1} (e^2-r^2) \right)  \]
with certain coefficients $c_k(j)$ independent of $d,e$.

Thirdly, by finding explicitly the value $\Ch_{2j-1}\big((-1) \times (2k-1)\big)$, we will determine the coefficients $c_k(j)$.

\medskip

The proofs of \cref{thm:odd,thm:even1,thm:even2}
are fully analogous to the proof of \cref{thm:main-odd-even} and we skip them, see \cref{sec:comments} for some additional comments.

\section{Characters on rectangular diagrams}
\label{sec:toolbox}

For any partition $\pi$ of $k$, Stanley's character formula for rectangular shapes \cite{Stanley2003/04} is given by
\begin{equation} 
    \label{eq:Stanley_rec}
    \Ch_\pi (p \times q)= (-1)^k\sum_{\substack{\sigma_1,\sigma_2 \in \Sym{k} \\ 
            \sigma_1 \sigma_2 = w_\pi}} 
    (-q)^{\kappa(\sigma_1)} p^{\kappa(\sigma_2)},
\end{equation}
where $\kappa(\sigma)$ is the number of cycles in $\sigma$ and $w_\pi$ is a fixed permutation of the cycle type $\pi$.

\begin{cor} 
    \label{cor:poly-degree}
For each partition $\pi$ the corresponding character
\[ \Ch_\pi (p \times q)\in\bZ[p,q]\] 
can be identified with a polynomial in the variables $p$ and $q$.
This polynomial is of degree $|\pi|+\ell(\pi)$ and fulfills the equality
\begin{equation}
\label{eq:transpose}    
\Ch_\pi ( p \times q) = (-1)^{|\pi|-\ell(\pi)} \Ch_\pi( q \times p ).
\end{equation}
\end{cor}
\begin{proof}
It is easy to show that if two polynomial functions from $\bQ[p,q]$ 
take equal values on each lattice point $(p,q)$ with $p,q\in\bN$ then they are equal as polynomials;
it follows that the polynomial given by the right-hand side of \eqref{eq:Stanley_rec}
is unique.

\medskip

We define the length $\| \sigma\|$ of a permutation $\sigma\in\Sym{k}$ as the minimal number of factors
necessary to write it as a product of transpositions. It is a classical result that
\[ \| \sigma \| = k - \kappa(\sigma). \] 
In this way the exponent on the right-hand side of \eqref{eq:Stanley_rec} is bounded from above by
\begin{multline*}
 \kappa(\sigma_1)+ \kappa(\sigma_2) = 2k - \|\sigma_1\| - \|\sigma_2 \|  \\ 
\leq 
2k - \| \sigma_1 \sigma_2 \|= k + \left( k - \| \sigma_1 \sigma_2 \| \right)
= |\pi|+\ell(\pi), 
\end{multline*}
as required.

\medskip

Equation \eqref{eq:transpose}  is a consequence of the general formula for the
character which corresponds to the transposed Young diagram
\[ \Ch_{\pi}(\lambda^T) = (-1)^{|\pi|-\ell(\pi)} \Ch_{\pi}(\lambda). \]
\end{proof}

\begin{cor}
\label{cor:minus-one}
For each integer $k\geq 1$
\begin{align*} 
\Ch_{k}\big((-1)  \times q\big) &= (-1)\ q (q+1) \cdots (q+k-1), \\
\Ch_{k}\big(p \times (-1) \big) &= (-1)^k\ p (p+1) \cdots (p+k-1).
\end{align*}
\end{cor}

\begin{proof}
We start with the following equality which holds in the symmetric group ring $\bZ[\Sym{k}]$
\begin{equation}  
\label{eq:factorization}
\sum_{\pi\in\Sym{k}} \pi = (1+J_1) (1+J_2) \cdots (1+J_k),   
\end{equation}
where
\[ J_i=(1,i)+\cdots+(i-1,i) \in \bZ[\Sym{k}] \]
denotes \emph{the Jucys--Murphy element}.
By investigating how the length of the permutations changes 
after multiplying by consecutive factors on the right-hand side of 
\eqref{eq:factorization} it follows from the Stanley formula \eqref{eq:Stanley_rec} that
\[ \Ch_{k}\big(p \times (-1) \big)  = 
(-1)^k \sum_{\sigma\in\Sym{k}} p^{\kappa(\sigma)} =
(-1)^k p (p+1) \cdots (p+k-1),
\]
as required.

\medskip

The other identity follows first one and \eqref{eq:transpose}.
\end{proof}

\section{Proof of \cref{thm:main-odd-even}}

We fix a positive integer $j$; note that the following notation depends on $j$ implicitly.
We denote by
\[\mathcal{P}(d)=\Ch_{2j-1} \big( (e-d) \times (e+d) \big)\]
the left-hand side of \eqref{eq:main-A}
\emph{viewed as a polynomial in the variable $d$} 
with the coefficients in the polynomial ring $\bQ[e]$. From the Stanley character formula \eqref{eq:Stanley_rec} and \cref{cor:poly-degree} it follows that the degree of $\mathcal{P}(d)$ is at most $2j$.

Equation \eqref{eq:transpose} implies that
\[ \Ch_{2j-1} \big( (e-d) \times (e+d) \big) =
 \Ch_{2j-1} \big( (e+d) \times (e-d) \big);\]
in other words the polynomial $\mathcal{P}(d)$ is even.

The linear space of even polynomials in the variable $d$ has a linear basis 
\[ 1, \quad d^2, \quad d^2 (d^2-1^2), \quad 
d^2 (d^2-1^2) (d^2-2^2), \quad \dots; \]
it follows that there exist polynomials $P_0, \ldots, P_{j} \in \bQ[e]$ with the property that
\begin{multline}
\label{eq:expansion}
\Ch_{2j-1} \big( (e-d) \times (e+d) \big)  \\ =
P_0(e) +
P_1(e)\ d^2 +
P_2(e)\ d^2 (d^2-1^2) +\cdots \\
=
\sum_{k=0}^j P_k(e)  \prod_{r=0}^{k-1} (d^2-r^2) .
\end{multline}
Additionally, from \cref{cor:poly-degree} 
it follows that the degree of the polynomial $P_k(e)$ is at most $2(j-k)$.

The parity of the total degree of each monomial 
on the right-hand side of the Stanley formula \eqref{eq:Stanley_rec} 
is the same as the parity of $|\pi|-\ell(\pi)$.
In our case $\pi=(2j-1)$ this parity is even;
it follows that
\begin{multline*} \Ch_{2j-1} \big( (e-d) \times (e+d) \big) =
\Ch_{2j-1} \big( (-e+d) \times (-e-d) \big)  
\\ = \Ch_{2j-1} \big( (-e-d) \times (-e+d) \big),
\end{multline*}
where the second equality is the consequence 
of the above observation that the polynomial $\mathcal{P}(d)$ is even. 
We proved in this way that $\mathcal{P}(d)$ is invariant under the involutive automorphism 
of the polynomial ring $\bQ[e]$ which is given by the substitution $e\mapsto -e$.
It follows that each coefficient $P_k(e)$ is an even polynomial in the variable $e$.

\medskip

\begin{lem}
\label{lem:induction}
For each $k\in\{0,\dots,j\}$
there exists some constant $c_k$ with the property that
\begin{equation} 
\label{eq:induction}
P_k(e) = c_k \prod_{r=k}^{j-1} (e^2-r^2) .
\end{equation}
\end{lem}
\begin{proof}
We will use induction over the variable $k$.
For the induction step let $k_0\in\{0,\dots,j\}$;
we assume that \eqref{eq:induction} holds true for each integer $k\in\{0,\dots,k_0-1\}$. 

Our strategy is to evaluate \eqref{eq:expansion} for $d=k_0$ and 
$e\in\{k_0-1,\dots,j-1\}$.
Each summand on the right-hand side  which corresponds to $k>k_0$ vanishes 
as it contains the factor $(d^2-r^2)$ for $r=k_0$.
On the other hand, 
each summand on the right-hand side which corresponds to $k<k_0$ vanishes because
either (a) $k_0=0$ and there are no such summands, or
(b) by the inductive hypothesis $P_k(e)$ contains the factor $(e^2-r^2)$ for $r=e$.
We proved in this way that for $e\in\{k_0-1,\dots,j-1\}$
\begin{equation}
\label{eq:pk0}
     \Ch_{2j-1} \big( (e-k_0) \times (e+k_0) \big)= 
P_{k_0}(e) \prod_{r=0}^{k_0-1} (k_0^2-r^2)  .
\end{equation}
In fact, in the special case when $k_0=0$ it is easy to check that
the above equality holds true for an \emph{arbitrary} choice of $e\in\bC$
and 
\begin{equation}
\label{eq:pk0prime}
     \Ch_{2j-1} \big( e \times e \big)= 
P_{0}(e).
\end{equation}
The formula \eqref{eq:pk0} has twofold consequences.

\medskip 

Firstly, in the special case when $e\in\{k_0,\dots,j-1\}$
the rectangular Young diagram $(e-k_0) \times (e+k_0)$
is well-defined and the defining formula \eqref{def:normalchar} 
can be used. Furthermore, this Young diagram 
does not contain any rim hooks of length $2j-1$; 
from the Murnaghan--Nakayama rule 
it follows that the left-hand side of \eqref{eq:pk0} is equal to zero;
as a consequence $P_{k_0}(e)=0$.

We proved in this way that $P_{k_0}$ is an even polynomial which has 
roots in $k_0,k_0+1,\dots,j-1$; it follows that the polynomial $P_{k_0}(e)$ is divisible
by the product
\[ \prod_{r=k_0}^{j-1} (e^2-r^2) . \]
Since the degree of $P_{k_0}$ is at most $2(j-k_0)$,
this determines the polynomial $P_{k_0}$ up to a scalar multiple
and shows that \eqref{eq:induction} holds true for $k:=k_0$.
This completes the proof of the inductive step of \cref{lem:induction}.
\end{proof}

As an extra bonus, for $k_0\geq 1$ 
the special case of \eqref{eq:pk0} and \eqref{eq:induction} 
for $e=k_0-1$ gives
(in order to keep the notation lightweight we write $k=k_0$)
\[ 
\Ch_{2j-1}\big( (-1) \times (2k-1) \big) =
c_{k} \prod_{r=k}^{j-1} \big((k-1)^2-r^2\big)
\prod_{r=0}^{k-1} (k^2-r^2). 
\]
The left-hand side can be evaluated thanks to \cref{cor:minus-one}
which gives an explicit product formula for the constant $c_k$ for $k\geq 1$.
Note that this argument cannot be applied in the special case $k=0$ 
because for $e=-1$ the information about the value of the polynomial 
$P_0(-1)$ is not linearly independent from the information about $P_0(1)=P_0(-1)$.

A combination of \eqref{eq:induction} and \eqref{eq:pk0prime}  gives
\[ \Ch_{2j-1} \big( e \times e \big) =   c_0 \prod_{r=0}^{j-1} (e^2-r^2).  \]
In order to evaluate the constant $c_0$ we need some additional piece of 
information about the polynomial on the left-hand side.
One possible approach is to evaluate its value for $e:=j$; 
in this special case the Murnaghan--Nakayama rule has only one summand therefore
value of the normalized character is given a product formula 
based on the hook-length formula. An alternative approach is based on 
calculating the leading coefficient $[e^{2j}] \Ch_{2j-1} \big( e \times e \big) $ 
based on the ideas of the asymptotic representation theory, 
see \cref{remark:asymptotic}.

Thanks to these explicit values of the constants $c_k$,
\cref{thm:main-odd-even} follows by a straightforward algebra
and its proof is now complete.

\medskip

\cref{thm:main-odd-even} gives a new proof of the following result. 

\begin{cor}
\label{cor:leading-Catalan}
For each integer $j\geq 1$ 
\[  \left[e^{2j} \right] \Ch_{2j-1}\big( e\times e \big) =
(-1)^{j-1} \operatorname{Cat}(j-1).\]
\end{cor}

\begin{remark}
\label{remark:asymptotic}
\cref{cor:leading-Catalan} is not new; 
in the following we only give a rough sketch of an alternative proof based on existing results.
The work of Biane (\cite[Theorem~1.3]{Biane1998} or 
\cite{Biane2003}) implies that 
\[ \lim_{e\to\infty} \frac{1}{e^{2j}} \Ch_{2j-1}\big(e \times e \big) = R_{2j}(\Box), \]
where $R_{2j}(\Box)$ denotes the \emph{free cumulant} of 
the one-box Young diagram $\Box=(1)$. More specifically, $R_{2j}(\Box)$ 
is the free cumulant of the \emph{Kerov transition measure} of $\Box$
which is equal to the Bernoulli measure
\[ \tfrac{1}{2} \left( \delta_{-1} + \delta_{1} \right).\]
Standard combinatorial tools of free probability
\cite{MingoSpeicher} give a closed formula 
for such a free cumulant in terms of Catalan numbers.
\end{remark}

\section{Comments about the proof of 
\cref{thm:odd,thm:even1,thm:even2}}
\label{sec:comments}

As we already mentioned, the proofs of 
\cref{thm:odd,thm:even1,thm:even2} are analogous to the proof \cref{thm:main-odd-even}. Below we revisit only some key places which require an adjustment.

\medskip

For example,
in order to prove \cref{thm:odd} we need to write
 \[ \Ch_{2j-1} \big( (e-d) \times (e+d) \big)  = 
 \sum_{k=0}^j P'_k(e)  \prod_{r=0}^{k-1} 
 \left(d^2-\left( r+\tfrac{1}{2}\right)^2\right) .
 \]
and then to show the following
analogue of \cref{lem:induction}: 
for each $k\in\{0,\dots,j\}$ there exists some constant $c_k'$ with the property that 
\[ P_k'(e)= c_k'  \prod_{r=k}^{j-1} \left(e^2-\left(r+\tfrac{1}{2}\right)^2 \right).\]
In order to achieve this goal, the strategy of the induction step is to fix $d=k_0+\tfrac{1}{2}$
and to consider $e\in\{ k_0-\tfrac{1}{2},\dots, j-\tfrac{1}{2} \}$.

The calculation of the constants $c'_k$ is particularly easy now because
both $c'_k$ as well as the constants $c_k$ from 
Equations \eqref{eq:induction} and \eqref{eq:expansion} 
coincide with the coefficient of a specific monomial in the Stanley
polynomial
\[ c'_k = \left[ d^{2k} e^{2(j-k)} \right] \Ch_{2j-1} \big( (e-d) \times (e+d) \big) = c_k
\]
hence they are equal.

\medskip

\cref{thm:even1,thm:even2} concern the character $\Ch_{2j}$ on an \emph{even} cycle.
In this case the corresponding polynomial 
\[\mathcal{P}(d)=\Ch_{2j} \big( (e-d) \times (e+d) \big)\]
is \emph{odd} and its degree is at most $2j+1$,
therefore we may write
\begin{align*}
\Ch_{2j} \big( (e-d) \times (e+d) \big) & = 
\sum_{k=0}^j P''_k(e)\  d \prod_{r=1}^{k} (d^2-r^2) \\
&= \sum_{k=0}^j P'''_k(e)\  d \prod_{r=1}^{k} \left( d^2-
\left( r-\tfrac{1}{2}\right) ^2 \right)
\end{align*}
for some even polynomials $P''_k(e), P'''_k(e)\in\bQ[e]$
which are of order at most $2(j-k)$,
where $k\in\{0,\dots,j\}$.

The proof of \cref{thm:even1} involves analysis of the polynomials $P''_k$ 
which is analogous to the one from the proof of \cref{thm:main-odd-even}; in particular
an analogue of \cref{lem:induction} says that
\[ P_k''(e)= c_k''  \prod_{r=k+1}^{j} \left(e^2-r^2 \right).\]
The proof of its inductive step is based on fixing $d=k_0+1$
and considering $e\in\{k_0,\dots,j\}$.
The values $e\in\{k_0+1,\dots,j\}$ are the positive roots of the even polynomial $P_k''$; the polynomial is therefore determined up to a multiplicative constant.
The special case $e=k_0$ allows to find explicitly the value of $c_k''$;
interestingly (opposite to the case in the proof of \cref{thm:main-odd-even})
the case $k=0$ does not require a separate proof.

\section{Integrality of the coefficients}

\begin{prop}
\label{prop:integer}
Let $d$ be an integer.
Each coefficient of the polynomial $\GA_d(j,n)$ (defined in \cref{cor:product}) 
is an integer. 
\end{prop}
\begin{proof}
We will show a stronger result that for each integer $k\geq 1$
\begin{equation}
\label{eq:my-dream}
\frac{1}{2^{k-1}} \frac{\prod_{r=0}^{k-1}(d^2-r^2)}{k! \, (2k-1)!!}
= \frac{2 d \cdot \prod_{r=-k+1}^{k-1} (d+r) }{(2k)!}
\end{equation}  
is an integer.
We will do it by proving that for each prime number $p$ the exponent by which
it contributes to the factorization of the numerator is at least its counterpart for the denominator. 
In the case when $p\neq 2$ these exponents are equal, respectively, to 
\begin{align} 
\label{eq:numerator}
 & \sum_{c\geq 1} \left( \big[ p^c \mid d \big] + 
\#\Big\{ i\in\{d-k+1,\dots,d+k-1\} : p^c \mid i \Big\} \right)  \\
\nonumber
  & = \sum_{c\geq 1} \left( \big[ p^c \mid d \big] - \big[ p^c \mid d+k \big] + 
\#\Big\{ i\in\{d-k+1,\dots,d+k\} : p^c \mid i \Big\} \right)
\intertext{and}
& \label{eq:denominator}
\sum_{c\geq 1}
\#\Big\{ i\in\{1,\dots,2k\} : p^c \mid i \Big\}.
\end{align}
Here and in the following we use the notation 
\[ [ \text{\emph{condition}} ] = 
\begin{cases} 1 & \text{if \emph{condition} holds true}, \\
0 & \text{otherwise}.
\end{cases}
\]
We will show that for each $c\geq 1$ the corresponding summand on the right-hand side of
\eqref{eq:numerator} is greater or equal to its counterpart in 
\eqref{eq:denominator}. 

\medskip

We start with the observation that in any collection of $p^c$ 
consecutive integers there is exactly one which is divisible by $p^c$; 
it follows that a collection of $2k$ consecutive integers contains at least $\left\lfloor \frac{2k}{p^c} \right\rfloor$ such numbers divisible by $p^c$.
As a consequence we get the following lower bound 
for the summand on the right-hand side of \eqref{eq:numerator}:
\begin{multline}
\label{eq:floorA}
\big[p^c \mid d \big] - \big[p^c \mid d +k \big] +
\#\Big\{ i\in\{d-k+1,\dots,d+k\} : p^c \mid i \Big\} \\
\geq [p^c \mid d ] - \big[p^c \mid d +k \big]  
+\left\lfloor \frac{2k}{p^c} \right\rfloor.
\end{multline}

Since \eqref{eq:numerator} can be alternatively written as
\begin{multline} 
\label{eq:numerator-B}
\sum_{c\geq 1} \left( \big[ p^c \mid d \big] + 
\#\Big\{ i\in\{d-k+1,\dots,d+k-1\} : p^c \mid i \Big\} \right)   \\
  = \sum_{c\geq 1} \left( \big[ p^c \mid d \big] - \big[ p^c \mid d-k \big] + 
\#\Big\{ i\in\{d-k,\dots,d+k-1\} : p^c \mid i \Big\} \right),
\end{multline}
an analogous reasoning to the one above gives the following alternative
lower bound for the summand on the left-hand side of \eqref{eq:numerator}
\begin{multline}
\label{eq:floorA2}
\big[p^c \mid d \big] - \big[p^c \mid d +k \big] +
\#\Big\{ i\in\{d-k+1,\dots,d+k\} : p^c \mid i \Big\}  \\
\geq [p^c \mid d ] - \big[p^c \mid d-k \big]  
+\left\lfloor \frac{2k}{p^c} \right\rfloor.
\end{multline}

\smallskip

On the other hand, for the corresponding summand in \eqref{eq:denominator} 
we get the following exact expression
\begin{equation}
\label{eq:floorB}
    \#\Big\{ i\in\{1,\dots,2k\} : p^c \mid i \Big\} = \left\lfloor \frac{2k}{p^c} \right\rfloor. 
\end{equation} 

\medskip

If one of the following two conditions holds true:
(a) the right-hand side of \eqref{eq:floorA} 
is greater or equal to the right-hand side of \eqref{eq:floorB}, or 
(b) the right-hand side of \eqref{eq:floorA2} 
is greater or equal to the right-hand side of \eqref{eq:floorB},
then the desired inequality holds true. 
The opposite case, 
i.e., when both $d-k$ and $d+k$ are divisible by $p^c$ 
and $d$ is not divisible by $p^c$, is clearly not possible since
$(d+k)+(d-k)=2d$.

\medskip

For the case when $p=2$ let us consider \emph{half} 
of the expression \eqref{eq:my-dream};
then \eqref{eq:numerator} and \eqref{eq:denominator} still provide the exponents
in the numerator and the denominator. 
The proof proceeds without any modifications until after Equation \eqref{eq:floorB}.
It might happen that none of the conditions (a) and (b) holds true;
if this is indeed the case then both $d-k$ and $d+k$ are divisible by $2^c$
(hence $d$ is divisible by $2^{c-1}$)
and $d$ is not divisible by $2^c$. 
There exists at most one index $c$ with this property.
It follows that the sum \eqref{eq:numerator} is at least 
\eqref{eq:denominator} \emph{less one}.
Since we considered \emph{half} of the expression \eqref{eq:my-dream}, 
this completes the proof.
\end{proof}

\begin{prop}
Let $d \in \{\pm \tfrac{1}{2}, \pm \tfrac{3}{2}, \dots\}$.
Each coefficient of the polynomial $\GB_d(j,n)$ (defined in \cref{cor:odd})
is an integer. 
\end{prop}

\begin{proof}
By the symmetry $\GB_d(j,n)=\GB_{-d}(j,n)$, we may suppose that $d$ is positive
and we write it as $d=d'+\tfrac{1}{2}$, where $d'$ is a non-negative integer. 
For each integer$k \in \{0, 1,2,\dots, d'\}$, 
\[
\frac{\prod_{r=0}^{k-1} \left( d^2- (r+\tfrac{1}{2})^2 \right)}{k! \, (2k-1)!!}
= \frac{2^k\prod_{r=-k+1}^{k}(d'+r)}{(2k)!} = 2^k \binom{d'+k}{2k}
\]
is clearly an integer.
Thus, we see that 
\begin{multline*}
\GB_{d'+\tfrac{1}{2}} (j,n)= \sum_{k=0}^{d'} (-1)^k 2^k \binom{d'+k}{2k}
j^{\downarrow k} (2j-1)^{\uparrow \uparrow k} \\
\times \prod_{r=k}^{d'-1} \left( n+ d'(d'+1) - r(r+1) \right)
\end{multline*}
has integer coefficient, as required.
\end{proof}

\begin{prop}
Let $d$ be an integer.
Each coefficient of the polynomial $\HA_d(j,n)$ (defined in \cref{cor:HA}) is an integer. 
\end{prop}

\begin{proof}
By the symmetry $\HA_d(j,n)=-\HA_{-d}(j,n)$, we may suppose that $d$ is a positive integer. 
For each integer $k \in \{1,2,\dots,d-1\}$,
\[
\frac{d\prod_{r=1}^k(d^2-r^2)}{k! \, (2k+1)!!}= \frac{2^k \prod_{r=-k}^{k} (d+r) }{(2k+1)!} = 2^k \binom{d+k}{2k+1}
\] 
is clearly an integer which completes the proof.
\end{proof}

\begin{prop}
Let $d$ be a half integer.
Each coefficient of the polynomial $\HB_d(j,n)$ (defined in \cref{cor:HB}) is an integer. 
\end{prop}

\begin{proof}
We may write $d=d'-\tfrac{1}{2}$, where $d'$ is an integer.
We will show a stronger result that for each integer $k\geq 1$
\[
\frac{1}{2^k} \frac{2d \prod_{r=1}^k \left( d^2 - (r-\tfrac{1}{2})^2 \right)}{k! \, (2k+1)!!} 
= \frac{(2d'-1) \prod_{r=-k}^{k-1} (d'+r)}{(2k+1)!}
\]
is an integer. The proof is analogous to the proof of \cref{prop:integer}.
\end{proof}

\section*{Acknowledgments} 
Sho Matsumoto was supported by JSPS KAKENHI, grant number 17K05281. 
Piotr \'{S}niady was supported by \emph{Narodowe Centrum Nauki}, grant number 
2017/26/A/ST1/00189.

\printbibliography

\end{document}